\documentclass[10pt]{article}

\usepackage{amsthm}
\usepackage{amsmath}
\usepackage{amsfonts}
\usepackage{amssymb}
\usepackage{mathrsfs}
\usepackage{empheq}
\usepackage{stmaryrd}
\usepackage{dsfont}
\usepackage{enumerate}
\usepackage{cancel}
\usepackage{comment}
\usepackage{pgfplots}
\usepackage{soul}
\usepackage{multirow}
\usepackage[applemac]{inputenc}
\usepackage{hyperref}
\usepackage{xcolor}

\definecolor{mygreen}{rgb}{0, 0.5, 0}

\newcommand{\1}{\mathds{1}}
\newcommand{\A}{\mathbb{A}}
\newcommand{\B}{\mathbb{B}}
\newcommand{\N}{\mathbb{N}}

\newcommand{\R}{\mathbb{R}}
\newcommand{\C}{\mathbb{C}}

\newcommand{\vsS}{\mathbb{S}}

\newcommand{\vH}{{\cal H}}

\newcommand{\vR}{{\cal R}}

\newcommand{\Dr}{\mathscr{D}}

\newcommand{\Lr}{\mathscr{L}}

\newcommand{\vphi}{\varphi}
\newcommand{\eps}{\varepsilon}
\newcommand{\dsp}{\displaystyle}
\newcommand{\ovl}{\overline}

\newcommand{\vlim}{\lim\limits}

\newcommand{\vint}{\int\limits}

\newcommand{\inj}{\hookrightarrow}

\newcommand{\weak}{\rightharpoonup}

\newcommand{\li}{\llbracket}
\newcommand{\ri}{\rrbracket}
\newcommand{\loc}{\mathrm{loc}}

\newcommand{\co}{\mathrm{c}}
\renewcommand{\d}{\mathrm{d}}

\newcommand{\dist}{\mathrm{dist}}

\newcommand{\M}{\mathrm{max}}

\newcommand{\CP}{C_\mathrm{P}}
\newcommand{\vi}{\mathrm{i}}
\newcommand{\w}{{\textsl w}}
\renewcommand{\le}{\leqslant}
\renewcommand{\ge}{\geqslant}
\renewcommand{\Re}{\mathrm{Re}}
\renewcommand{\Im}{\mathrm{Im}}

\newcommand{\p}{\prime}

\newcommand{\eqdef}{\stackrel{\mathrm{def}}{=}}
\DeclareMathOperator{\supp}{supp}

\numberwithin{equation}{section}

\newtheorem{thm}{Theorem}[section]

\newtheorem{lem}[thm]{Lemma}
\newtheorem{sym}[thm]{Symmetry Property}

\newtheorem{vthm}{Theorem}[subsection]
\newtheorem{vlem}[vthm]{Lemma}

\theoremstyle{definition}
\newtheorem{rmk}[thm]{Remark}
\newtheorem{defi}[thm]{Definition}

\newenvironment{proof*}{\noindent{\bf Proof.}}{\qed}
\newenvironment{vproof}[1]{\noindent{\bf Proof #1}}{\qed}

\title{\huge Solutions with expanding compact support of saturated Schrödinger equations: self-similar solutions}
\author{\sc Pascal Bégout$^*$ and Jes\'us Ildefonso D\'{\i}az$^\dagger$}
\date{}

\begin{document}

\maketitle

\begin{center}
\begin{tabular}{ll}
\hspace*{-.28cm}$^*$ Toulouse School of Economics	&	\hspace*{-.25cm}$^\dagger$ Instituto de Matem\'atica Interdisciplinar	\\
Universit\'e Toulouse Capitole 						&	Universidad Complutense de Madrid			\\
Institut de Math\'ematiques de Toulouse 				&	Plaza de las Ciencias, 3						\\
1, Esplanade de l' Universit\'e 						&	28040 Madrid, SPAIN						\\
31080 Toulouse Cedex 6, FRANCE					& 	\\
{\footnotesize E-mail\:: \href{mailto:Pascal.Begout@math.cnrs.fr}{\texttt{Pascal.Begout@math.cnrs.fr}}}
&
{\footnotesize E-mail\:: \href{mailto:jidiaz@ucm.es}{\texttt{jidiaz@ucm.es}}}
\end{tabular}
\end{center}

\begin{abstract}
We prove the existence of some solutions $u(t,x)$  of the Schrödinger equation with a saturation nonlinear term $(u/|u|)$ having compact support, for each $t>0,$ expanding with $t$ with a growth law of the type $C\sqrt{t}.$ The primary tool is considering the self-similar solution of the associated equation.
\end{abstract}

{\let\thefootnote\relax\footnotetext{Pascal Bégout acknowledges funding from ANR under grant ANR-17-EUR-0010 (Investissements d’Avenir program)}}
{\let\thefootnote\relax\footnotetext{The research of J.I.~D\'{\i}az was partially supported by the project PID2023-146754NB-I00 funded by MCIU/AEI/10.13039/501100011033 and FEDER, EU. MCIU/AEI/10.13039/-501100011033/FEDER, EU.}}
{\let\thefootnote\relax\footnotetext{$^*$\href{https://orcid.org/0000-0002-9172-3057}{https://orcid.org/0000-0002-9172-3057}}}
{\let\thefootnote\relax\footnotetext{$^\dagger$\href{https://orcid.org/0000-0003-1730-9509}{https://orcid.org/0000-0003-1730-9509}}}
{\let\thefootnote\relax}
{\let\thefootnote\relax}

\tableofcontents

\baselineskip .6cm

\section{Introduction}
\label{introduction}

The existence of compactly support solutions to Schrödinger equation was a constant subject of research since Schrödinger postulated the existence of such equation in 1925 and published it in 1926. For the case of the linear equation it seems that it was Sir Nevill Francis Mott (1905-1996), who would
later win the Nobel Prize in 1977, proposed the study of the infinite well potential in his 1930 book \cite{Mott}. This was a generalization of the finite well potential proposed, in 1928, by George Gamow \cite{Gamow} when finding the tunnel effect by first time in the literature. Solutions of the linear Schrödinger equation with an infinite well potential have compact support (the compact set of $\R^N$ where the potential is finite) but the mathematical study of this problem presents some ambiguities (\cite{MR3421910}) which disappear when such a discontinuous potential is replaced by strongly singular potentials of the Pöschl--Teller type (\cite{PoschTeller}, \cite{MR3421910}, \cite{MR3721865}).
\medskip \\
This study of the support of solutions of nonlinear Schrödinger equations was also considered by many authors but with negative results when the nonlinear term is Lipschitz continuous (see, e.g., the presentation made by J. Bourgain in \cite{MR1443322}). These authors made completely new contributions in the subject by showing that solutions with compact support do exist when the nonlinear term is not Lipschitz continuous but of the form
\begin{equation}
\label{EqSingularSaturation}
\vi\dfrac{\partial u}{\partial t}+\Delta u=a|u|^{-(1-m)}u+f(t,x),
\end{equation}
for some $m\in(0,1)$ and for a suitable complex coefficient $a.$ This equation is associated to the consideration of non-Kerr law optical Schrödinger equation arising, for instance, in nonlinear optical media. This type of equations also arises in Quantum Mechanics and Hydrodynamics. When
searching for ``solitary wave solutions" of the form $u(t,x)=\psi(x)e^{\vi bt}$ (when $f(t,x)=e^{\vi bt}F(x))$ then the complex function $u$ satisfies a stationary nonlinear equation which leads to solutions with compact support once we assume that $F(x)$ has compact support.
\medskip \\
The above problem was extended to the case of saturated nonlinear terms $(m=0)$ in the recent paper \cite{MR4848642} proving that ``solitary wave solutions" $u(t,x)=\psi(x)e^{\vi bt}$ have compact support even if $F(x)$ does not have compact support but is small enough outside of some compact subset of $\R^N.$ From the qualitative point of view, the above type of solutions with compact support (for the mentioned linear and nonlinear cases) concern some special type of solutions: ``solitary wave solutions" of the form $u(t,x)=\psi(x)e^{\vi bt}$ which implies that support of $u(t)$ does not move, for any $t>0,$ since $\supp u(t)=\supp\psi.$
\medskip \\
A different point of view was followed by the authors in the paper \cite{MR3193996} where the existence of a self-similar solution of the form $u(t,x)=t^{\frac{p}2}\vphi\left(\frac{x}{\sqrt{t}}\right)$ was proved for the equations of the type \eqref{EqSingularSaturation} with $m\in (0,1)$ once we assume that $f(t,x)=t^{\frac{p-2}2}F\left(\frac{x}{\sqrt{t}}\right):$ it was proved in that paper that if $\supp F$ is compact then the solution profile $\vphi$ is also compact. As it was detailed later, for this type of solutions their support $\supp u(t)$ expands with time $t>0,$ with a sublinear growth of the type $C\sqrt{t}.$
\medskip \\
The main objective of the present paper is to extend the results of \cite{MR3193996} to the saturated case $(m=0)$ by showing that the corresponding solution has expanding support $\supp u(t)$ that expands with time $t>0,$ with a sublinear growth of type $C\sqrt{t}$ even if the profile $F$ of the data $f(t,x)=t^{\frac{p-2}2}F\left(\frac{x}{\sqrt{t}}\right) $ is not compactly supported but, as in \cite{MR4848642}, is sufficiently small outside a compact subset of $\R^N.$ One of the consequences of such general assumption on $f(t,x)$ is that we can extend the property of solutions with compact support when we couple the Schrödinger equation with some other phenomena (as for instance the existence of some magnetic fields: see Section~9 of \cite{MR4848642}). Here we are interested in finding self-similar solutions with compact support in the space variable of the following Schrödinger equation with saturated nonlinearity,
\begin{gather}
\label{nlse}
\begin{cases}
\vi\dfrac{\partial u}{\partial t}+\Delta u=a\,U+f(t,x), \, (t,x)\in(0,\infty)\times\R^N,	\medskip \\
U=\dfrac{u}{|u|}, \text{ a.e.\:in } \big\{(t,x)\in(0,\infty)\times\R^N;u(t,x)\neq0\big\},
\end{cases}
\end{gather}
where $a\in\C.$ For this, it is enough to study the equation satisfied by the profile $\vphi$ of $u$, that is
\begin{gather}
\label{U}
\begin{cases}
-\Delta\vphi + a\,\Phi - \dfrac{\vi p}2\vphi + \dfrac{\vi}2x.\nabla \vphi = -F, \text{ in } \Dr^\p(\R^N), \medskip \\
\Phi=\dfrac\vphi{|\vphi|}, \text{ a.e.\:in } \big\{x\in\R^N;\vphi(x)\neq0\big\},
\end{cases}
\end{gather}
where $p\in\C$ with $\Re(p)=2,$ $\vphi=u(1),$ and $F=f(1).$ We will maintain the
notation and several common arguments with our previous papers \cite{MR3193996} and \cite{MR4848642} but new results will be given improving both papers. As in \cite{MR3193996}, it is useful to introduce a change of unknowns which brings us back to the search for solutions to the problem
\begin{gather}
\label{g}
\begin{cases}
-\Delta g + a\,G - \vi\dfrac{N+2p}4g - \dfrac1{16}|x|^2g = -Fe^{-\vi\frac{|x|^2}8}, \text{ in } \Dr^\p(\R^N),	\medskip \\
G=\dfrac{g}{|g|}, \text{ a.e.\:in } \big\{x\in\R^N;g(x)\neq0\big\}.
\end{cases}
\end{gather}
So in this paper, we study the following which is more general equation than \eqref{g},
\begin{gather}
\label{nls}
\begin{cases}
-\Delta u+a\,U+b\,u+V\,u=F, \text{ in } H^{-1}(\Omega)+L^\infty(\Omega),	\medskip \\
U=\dfrac{u}{|u|}, \text{ a.e.\:in } \omega=\big\{x\in\Omega;u(x)\neq0\big\},
\end{cases}
\end{gather}
where $(a,b)\in\C^2,$ and $\Omega$ is a subset of $\R^N$ whose boundary is $\Gamma,$ with homogeneous Dirichlet boundary condition
\begin{gather}
\label{dir}
 u_{|\Gamma}=0,
\end{gather}
or with homogeneous Neumann boundary condition
\begin{gather}
\label{neu}
 \dfrac{\partial u}{\partial\nu}_{|\Gamma}=0.
\end{gather}
The compactness of the support of solutions will be obtained by some improvements of the energy methods presented in the monograph \cite{MR2002i35001} (see also the extension to some variational inequalities made in \cite{MR2466410}). We mention that the method such as it was developed in the above mentioned references is only well adapted to its application to complex problems of Ginzburg -Landau type \cite{MR3208711} in which the time derivative of the unknown contains a real part (situation which is not valid for the Schrödinger equation).
\medskip \\
The organization of this paper is as it follows: Section~\ref{secsss} is devoted to the structure of self-similar solutions and to the presentation of the main result of this paper (Theorem~\ref{thmmain} below). Details on the notion of solutions, the results on the existence and uniqueness of solutions are collected in Section~\ref{seceu}. A set of auxiliary results preparing the application of an energy method leading to the compactness of the support of the solution, as well as the proof of the results stated in the previous sections are presented in Section~\ref{secset}. Finally, an Appendix is devoted to the proof of the additional regularity obtained from the structure of self-similar solutions.
\medskip \\
As indicated before, this paper extends some previous papers by the authors (\cite{MR3193996} and a part of \cite{MR3315701}) to the case $m=0.$ Nevertheless, since the applied techniques are of some different type, they do not allow to conclude some previous results in their complete generality. Furthermore, despite the fact that \cite{MR4848642} also concerns equation \eqref{nls}, we point out that the assumptions and results differ so that they cannot be employed to construct self-similar solutions with compact support in space. For instance, Theorem~\ref{thmexi1} vs \cite[Theorem~2.6]{MR4848642}: \cite[Theorem~2.6]{MR4848642} in less restrictive in terms of $V,$ and Theorem~\ref{thmexi1}, only considers the Dirichlet condition and $|\Omega|<\infty.$ But Theorem~\ref{thmexi1} is more general in terms of $(a,b)$ since $(a,b)\in\C\times\B$ while in \cite[Theorem~2.6]{MR4848642}, $(a,b)\in\A^2$ satisfy some additional conditions, and $\A\subsetneq\B.$ Theorem~\ref{thmexi2} vs \cite[Theorem~2.6]{MR4848642}: Theorem~\ref{thmexi2} is more restrictive in terms of $(a,b)$ but it allows $V$ to be a complex-valued function with no sign restriction about $\Re (V),$ while in \cite[Theorem~2.6]{MR4848642}, $V$ is a nonnegative real-valued functions. It is essential to allow to choose $V$ with a negative real part to consider self-similar solutions.
\medskip \\
Here is a list of symbols we will use in this paper: for a complex number $z,$ we denote by $\ovl z,$ $\Re(z)$ and $\Im(z),$ its conjugate, real and imaginary part, respectively, and $\vi^2=-1.$ $\N_0=\N\cup\{0\}.$ For $p\in[1,\infty],$ $p^\prime$ is the conjugate of $p$ defined by $\frac{1}{p}+\frac{1}{p^\prime}=1.$ Unless specified, all functions are complex-valued and all the vector spaces are considered over the field $\R.$ For a Banach space $X,$ we denote by $X^\star\eqdef\Lr(X;\R)$ its topological dual and by $\langle\: . \; , \: . \:\rangle_{X^\star,X}$ the $X^\star-X$ duality product. By convention, $W^{0,q}(\R^N)=L^q(\R^N),$ for any $0<q<\infty.$ For positive parameters $a_1,\ldots,a_n,$ we shall write $C(a_1,\ldots,a_n)$ to indicate that $C$ is a positive constant which depends only and continuously on $a_1,\ldots,a_n.$ Finally, if $A$ is a subset of $\R^N$ then $A^\co$ denotes its complement, and $A\setminus B=A\cap B^\co.$
\medskip \\
Let us recall that if $X$ and $Y$ are two Banach spaces\footnote{Actually, locally convex topological vector spaces is enough which allows to consider $X=\Dr(\Omega).$} such that $X\inj Y$ with dense embedding then $Y^\star\inj X^\star,$ and for any $F\in Y^\star$ and $u\in X,$ $\langle F,u\rangle_{X^\star,X}=\langle F,u\rangle_{Y^\star,Y}.$ By the Riesz representation Theorem, we have for any $p\in[1,\infty),$ $F\in L^{p^\p}(\Omega)$ and $u\in L^p(\Omega),$ $\langle F,u\rangle_{L^{p^\p}(\Omega),L^p(\Omega)}=\Re\int_\Omega F(x)\ovl{u(x)}\d x.$ In particular, this implies that we shall always identify $L^2(\Omega)$ with its topological dual. In addition, if $A_1$ and $A_2$ are two Banach spaces such that $A_1,A_2\subset\vH$ for some Hausdorff topological vector space $\vH,$ and if $A_1\cap A_2$ is dense in both $A_1$ and $A_2$ then $A_1\cap A_2$ and $A_1+A_2$ are Banach spaces, and $\big(A_1\cap A_2\big)^\star=A_1^\star+A_2^\star.$ This justifies the identity~\eqref{defsol11} below. For more details, see Trèves~\cite{MR2296978}, Bergh and Löfström~\cite{MR0482275}, and~\cite{MR4521439}.

\section{Self-similar solutions}
\label{secsss}

\noindent
Let us recall that the notion of self-similar solutions relies on the transformation $\lambda\longmapsto(u_\lambda,U^\lambda),$ where for $\lambda>0,$ $p\in\C,$ $u\in L^1_\loc\big((0,\infty)\times\R^N\big)$ and $U$ a saturated section associated to $u$ (Definition~\ref{defsat} below),
\begin{gather}
\label{ula}
u_\lambda(t,x)=\lambda^{-p}u(\lambda^2t,\lambda x),	\\
\label{Ula}
U^\lambda(t,x)=\lambda^{-(p-2)}U(\lambda^2t,\lambda x),
\end{gather}
for a.e.\,$(t,x)\in(0,\infty)\times\R^N.$ We also recall that $\lambda^p\eqdef e^{p\ln\lambda}$ and $|\lambda^p|=\lambda^{\Re(p)}.$ If $\Re(p)=2$ then a straightforward calculation shows that if $(u,U)$ is a solution to \eqref{nlse} with $f=0,$ then so is $(u_\lambda,U^\lambda),$ for any $\lambda>0.$ In particular, $U^\lambda$ is a saturated section associated to $u_\lambda.$ To keep this property when $f\neq0,$ with $f\in L^1_\loc\big((0,\infty)\times\R^N\big),$ we assume that $f$ satisfies
\begin{gather}
\label{fla}
\forall\lambda>0, \; f^\lambda=f,
\end{gather}
or equivalently,
\begin{gather}
\label{eqprof}
f(t,x)=t^{\frac{p-2}2}F\left(\frac{x}{\sqrt t}\right),
\end{gather}
for a.e.\,$(t,x)\in(0,\infty)\times\R^N,$ where $F=f(1).$ To have functions $f$ satisfying \eqref{fla}, it is sufficient for any given function $F\in L^1_\loc(\R^N)$ to define $f$ by \eqref{eqprof}. Furthermore, we easily check that $(u,U)$ satisfies the invariance property
\begin{gather*}
\forall\lambda>0, \; (u_\lambda,U^\lambda)=(u,U),
\end{gather*}
if, and only if,
\begin{gather}
\label{eqprou}
u(t,x)=t^\frac{p}2\vphi\left(\frac{x}{\sqrt t}\right),	\\
\label{eqproU}
U(t,x)=t^\frac{p-2}2\Phi\left(\frac{x}{\sqrt t}\right),
\end{gather}
for a.e.\,$(t,x)\in(0,\infty)\times\R^N,$ where $(\vphi,\Phi)=(u(1),U(1)).$ This remarkable invariance property leads to the well-known definition of self-similar solution.

\begin{defi}
\label{defsat}
Let $\theta\subseteq\R^N$ be an open subset and let $u\in L^1_\loc(\theta).$ A function $U\in L^\infty(\theta)$ is said to be a \textit{saturated section} associated to $u$ if $\|U\|_{L^\infty(\theta)}\le1$ and $U=\frac{u}{|u|},$ almost everywhere in $\omega\eqdef\big\{y\in\theta;u(y)\neq0\big\}.$
\end{defi}

\begin{defi}
\label{defselsim}
Let $f\in C\big((0,\infty);L^2(\R^N)\big)$ satisfy~\eqref{fla} and let $p\in\C$ be such that $\Re(p)=2.$ A solution $(u,U)$ to~\eqref{nlse} is said to be \textit{self-similar} if $u\in C\big((0,\infty);L^2(\R^N)\big),$ $U$ is a saturated section associated to $u$ and if for any $\lambda>0,$ $(u_\lambda,U^\lambda)=(u,U),$ where $u_\lambda$ and $U^\lambda$ are defined by~\eqref{ula} and \eqref{Ula}, respectively. In this cases, $u(1)$ is called the \textit{profile} of $u$ and is denoted by $\vphi.$
\end{defi}

\noindent
It follows from~\eqref{nlse}, \eqref{eqprou} and~\eqref{eqproU} that the profile $\vphi$ of $u$ and $\Phi$ satisfy \eqref{U}. In particular, $\Phi$ is a saturated section associated to $\vphi.$ Conversely, if $(\vphi,\Phi)\in L^2(\R^N)\times L^\infty(\R^N)$ satisfies~\eqref{U} with $\|\Phi\|_{L^\infty(\R^N)}\le1,$ then the functions $u$ and $U$ defined by \eqref{eqprou} and \eqref{eqproU}, respectively, belong to $C\big((0,\infty);L^2(\R^N)\big)$ (Lemma~\ref{lemA}) and $L^\infty((0,\infty)\times\R^N),$ respectively, $U$ is a saturated section associated to $u$ and $u$ is a self-similar solution to~\eqref{nlse}, where $f$ is defined by \eqref{eqprof} and satisfies \eqref{fla}. A priori estimates on $\vphi$ are not easy to obtain due to the term $x.\nabla\vphi.$ Thus in the literature, this problem is circumvented using the bijective transformation
\begin{gather}
\label{gU}
g(x)=\vphi(x)e^{-\vi\frac{|x|^2}8}, \; \text{for a.e.\,} x\in\R^N.
\end{gather}
The saturated section $\Phi$ associated to $\vphi$ then becomes
\begin{gather}
\label{gP}
G(x)=\Phi(x)e^{-\vi\frac{|x|^2}8}, \; \text{for a.e.\,} x\in\R^N.
\end{gather}
It follows that for any $p\in\C$ and $\vphi\in L^2(\R^N),$ whose a saturated section associated to $\vphi$ is $\Phi,$ $(\vphi,\Phi)$ is a solution to~\eqref{U} if, and only if, $(g,G)\in L^2(\R^N)\times L^\infty(\R^N)$ is a solution to \eqref{g} and $G$ is a saturated section associated to $g.$ The study of \eqref{g} is then more convenient than of \eqref{U}, and is related to Theorem~\ref{thmexi2}. Let,
\begin{gather}
\label{A}
\A=\C\setminus\big\{z\in\C; \Re(z)\le0 \text{ and } \Im(z)=0\big\}.
\end{gather}
The main result of this paper is the following.

\begin{thm}
\label{thmmain}
Assume that $a\in\A$ is such that $\Im(a)\le0.$ Let $p\in\C$ be such that $\Re(p)=2,$ let $f\in C\big((0,\infty);L^2(\R^N)\big)$ satisfy~\eqref{fla} and set $F=f(1).$ Assume also that $F_{|K^\co}\in L^\infty(K^\co),$ for some compact subset $K$ of $\R^N.$
\begin{enumerate}
\item
\label{thmmain1}
\textbf{Existence.}
For any $R>0$ such that $K\subset\ovl B(0,R)$ and any $\eps>0,$ there exist $M=M(|a|,|\Im(p)|,R,N)$ and $\delta=\delta(|a|,|\Im(p)|,R,\eps,N)$ satisfying the following property. If $\|F\|_{L^2(\R^N)}\le\delta$ and $\|F\|_{L^\infty(K^\co)}\le\frac1M$ then there exists a self-similar solution $(u,U)$ to~\eqref{nlse} such that
\begin{gather}
\label{thmmain11}
u\in C\big((0,\infty);H^2(\R^N)\big)\cap C^1\big((0,\infty);H^1(\R^N)\big)\cap C^2\big((0,\infty);L^2(\R^N)\big)
\end{gather}
and for any $t>0,$ $\supp u(t)$ is compact. In addition, the profile $\vphi$ of $u$ satisfies that $\supp\vphi\subset K(\eps)\subset\ovl B(0,R+\eps),$ where
\begin{gather*}
K(\eps)=\Big\{x\in\R^N;\; \dist(x,K)\le\eps\Big\},
\end{gather*}
which is compact.
\item
\label{thmmain2}
\textbf{Uniqueness.}
Let $(u,U)$ and $(v,V)$ be two self-similar solutions to \eqref{nlse} with profiles $\vphi$ and $\phi,$ respectively, and with $\supp\vphi\cup\supp\phi\subset B(0,r),$ for some $r>0.$ Assume that one of the two following conditions is satisfied.
\begin{enumerate}
\item
\label{thmmain21}
$\Re(a)=0.$
\item
\label{thmmain22}
$\Re(a)>0$ and $r^2\le8\Im(p)+4\frac{|\Im(a)|}{\Re(a)}(N+4).$
\end{enumerate}
Then for any $t>0,$ $u(t)=v(t).$ As a consequence, $U=V$ almost everywhere in $(0,\infty)\times\R^N.$
\end{enumerate}
\end{thm}

\noindent
We postpone the proof of Theorem~\ref{thmmain} to Subsection~\ref{subsecproof}.

\begin{rmk}
\label{rmkuni}
It is obvious from \eqref{nlse} that the uniqueness of the solution $u$ implies the uniqueness of the saturated section $U.$
\end{rmk}

\begin{rmk}
\label{rmkthmmain}
In \cite{MR3193996}, self-similar solutions are studied with the nonlinearity $|u|^{-(1-m)}u,$ where $0<m<1.$ It is shown that a self-similar solution cannot be continuous at $t=0$ in a reasonable way. This remains true in our case (which corresponds to $m=0).$ Below, we give some details. Let $p\in\C$ be such that $\Re(p)=2$ and let $u$ be a self-similar solution to~\eqref{nlse} with profile $\vphi.$
\begin{enumerate}
\item
\label{rmkthmmain1}
Let us define the transformation $T_\lambda:v\longmapsto v_\lambda,$ for any $v\in L^1_\loc(\R^N),$ when $\lambda>0:$ $T_\lambda(v)(\,.\,)=\lambda^{-p}v(\lambda\:.\:).$ The functions which satisfy this invariance property cannot be $L^q$-functions in the sense that we have
\begin{gather*}
\Lambda_q\eqdef\big\{v\in L^q(\R^N);\forall\lambda>0, \; T_\lambda(v)=v\big\}=\big\{0\big\},
\end{gather*}
for any $q\in(0,\infty].$ Indeed, if for some $q\in(0,\infty],$ $v\in\Lambda_q$ then a straighforward calculation gives that
\begin{gather*}
\forall\lambda>0, \; \|v\|_{L^q(\R^N)}=\lambda^{2+\frac{N}q}\|v\|_{L^q(\R^N)}.
\end{gather*}
Therefore, $v=0.$ It follows that if $u(0)\in L^q(\R^N),$ for some $0<q\le\infty,$ then $u(0)\in\Lambda_q$ and so necessarily $u(0)=0.$
\item
\label{rmkthmmain2}
It follows from above that if $u\in C\big([0,\infty);\Dr^\p(\R^N)\big)$ is a self-similar solution to~\eqref{nlse} with $u(0)\neq0$ then for any $0<q\le\infty,$ $u\notin C\big([0,\infty);L^q(\R^N)\big).$ On the other hand, if for some $0<q\le\infty,$ $u\in C\big((0,\infty);L^q(\R^N)\big)$ then $\vphi\in L^q(\R^N)$ and it follows from \eqref{eqprou} that
\begin{gather}
\label{rmkthmmain21}
\forall t>0, \; \|u(t)\|_{L^q(\R^N)}=t^{1+\frac{N}{2q}}\|\vphi\|_{L^q(\R^N)},
\end{gather}
and so $\vlim_{t\searrow0}\|u(t)\|_{L^q(\R^N)}=0.$ Actually, if $m\in\{0,1,2\},$ $0<q\le\infty$ and $\vphi\in W^{m,q}(\R^N)$ then by \eqref{eqprou}, $u(t)\in W^{m,q}(\R^N),$ for any $t>0,$ and
\begin{gather}
\label{rmkthmmain22}
\|\nabla u(t)\|_{L^q(\R^N)}=t^{\frac12+\frac{N}{2q}}\|\nabla\vphi\|_{L^q(\R^N)},	\\
\label{rmkthmmain23}
\|\partial^2_{jk}u(t)\|_{L^q(\R^N)}=t^\frac{N}{2q}\|\partial^2_{jk}\vphi\|_{L^q(\R^N)},
\end{gather}
for any $t>0$ and $(j,k)\in\li1,N\ri^2,$ so that $\vlim_{t\searrow0}\|u(t)\|_{W^{m,q}(\R^N)}=0$ $(q<\infty,$ if $m=2).$
\item
\label{rmkthmmain3}
If $f=0,$ $a\in\R$ and $\vphi$ has compact support then for any $t\in\R,$ $u(t)=0.$ Indeed, if $g$ is defined by \eqref{gU} then $g\in L^2(\R^N)$ and by \eqref{g}, $\Delta g\in L^2_\loc(\R^N).$ By interior elliptic regularity, $g\in H^2_\loc(\R^N)$ (Cazenave~\cite[Proposition~4.1.2]{caz-sle}). Then $\vphi\in H^2_\loc(\R^N)$ and since $\supp\vphi$ is compact, we finally have $\vphi\in H^2(\R^N).$ It follows from Lemma~\ref{lemA} below that $u$ satisfies the regularity \eqref{thmmain11}. We are then allowed to take the $X^\star-X$ duality product of \eqref{nlse} with $\vi u,$ where $X=H^1(\R^N)\cap L^1(\R^N),$ to obtain that $\frac{\d}{\d t}\|u(t)\|_{L^2(\R^N)}^2=0,$ for any $t>0.$ With help of \eqref{rmkthmmain21}, we then deduce that
\begin{gather*}
\forall t>0, \; \|\vphi\|_{L^2(\R^N)}=\|u(t)\|_{L^2(\R^N)}=t^{1+\frac{N}4}\|\vphi\|_{L^2(\R^N)}.
\end{gather*}
Then $\vphi=0,$ from which the result follows.
\item
\label{rmkthmmain4}
Assume that $u(0)\neq0.$ From the structure of the self-similar solution $u$ we easily deduce that for any $t>0,$
\begin{gather*}
\supp u(t)=\sqrt t\,\supp\vphi.
\end{gather*}
Letting $t\searrow0,$ we could conclude that $\supp u(0)=\emptyset$ and then $u(0)=0.$ But as seen above, $u$ is not continuous at $t=0$ in any reasonable way and we cannot infer that $u(0)=0.$ Estimates on the expansion of the support of the type $C\sqrt{t}$ where proved, for the first time, for parabolic variational inequalities, in the paper H.~Brezis and A.~Friedman~\cite{MR390501}.
\end{enumerate}
\end{rmk}

\begin{rmk}
\label{rmkild2}
Let $0<m<1,$ let $a\in\A$ be such that $\Im(a)\le0,$ let $p\in\C$ be such that $\Re(p)=\frac2{1-m}$ and let $f_1,\ldots,f_d\in C\big((0,\infty);L^2(\R^N)\big)$ satisfying~\eqref{fla}. Assume further that for any $j\in\li1,d\ri,$ $K_j\eqdef\supp f_j(1)$ is compact, $\|f_j(1)\|_{L^2(\R^N)}$ is small enough and $K_j\cap K_\ell=\emptyset,$ for any $j\neq\ell.$ It follows from \cite[Theorem~1.2]{MR3193996} that for any $j\in\li1,d\ri,$ there exists a self-similar solution $u_j$ to
\begin{gather*}
\vi\dfrac{\partial u_j}{\partial t}+\Delta u_j=a|u_j|^{-(1-m)}u_j+f_j(t,x), \, (t,x)\in(0,\infty)\times\R^N,
\end{gather*}
such that $\supp u_j(1)$ is compact. Due to the smallness of the $d$ norms $\|f_j(1)\|_{L^2(\R^N)},$ we also have that for any $j\neq\ell,$ $\supp u_j(1)\cap\supp u_\ell(1)=\emptyset.$ Set,
\begin{gather*}
u=\sum_{j=1}^du_j \; \text{ and } \; f=\sum_{j=1}^df_j.
\end{gather*}
From the structure of the self-similar solutions and since the support of the $d$ functions $u_j(1)$ are disjoints, we conclude that the supports of the $d$ functions $u_j$ remain disjoints at least during some suitable period of time $(0,T),$ for some $T>1.$ It follows that $u$ is a self-simlar solution to
\begin{gather*}
\vi\dfrac{\partial u}{\partial t}+\Delta u=a|u|^{-(1-m)}u+f(t,x), \, (t,x)\in(0,T)\times\R^N,
\end{gather*}
 although this equation is not linear. If $m=0,$ then the above arguments do not work since we do not necessarily have that the saturated section $U$ associated to $u$ satisfies $U=0$ when $u=0.$ Nevertheless, we may still generate self-similar solutions of the evolution Schrödinger equation~\eqref{nlse}, on a finite time interval $(0,T),$ with $T>1,$ having support with more than one connected component. Indeed, it is sufficient to work with one function $f$ where the compactness of $f(1)$ and the smallness of $\|f(1)\|_{L^2(\R^N)}$ are replaced by the following assumptions: there exist $d$ compact connected subsets $K_j$ such that $\|f(1)\|_{L^\infty(K^\co)}$ is small enough, where $K=\cup_{j=1}^dK_j,$ and such that for any $j\neq\ell,$ $K_j\cap K_\ell=\emptyset.$ We conclude with help of Theorem~\ref{thmmain}.
\end{rmk}

\begin{rmk}
\label{rmkild3}
It is useful to rewrite the evolution Schrödinger equation in terms of real components of solutions and data 
\begin{gather*}
\begin{cases}
u=u_R+\vi u_I,		&	f=f_R+\vi f_I,	\medskip \\
a=a_R+\vi a_I.		& 
\end{cases}
\end{gather*}
Then, the sign of the components of the coefficient $a$ is especially crucial for understanding the different nature of the coupled system. Theorem~\ref{thmmain} holds, for instance, if $a=\lambda-\vi\mu $ with $\lambda,\mu >0$ (in the pure elliptic system, the case of $\mu <0$ is also allowed: see our paper in~\cite{MR4848642}), and then we arrive to the coupled system 
\begin{gather*}
\begin{cases}
\dfrac{\partial u_I}{\partial t}-\Delta u_R+\dfrac{\lambda u_R+\mu u_I}{\sqrt{u_R^2+u_I^2}}=-f_R,	\medskip \\
-\dfrac{\partial u_R}{\partial t}-\Delta u_I+\dfrac{\lambda u_I-\mu u_R}{\sqrt{u_R^2+u_I^2}}=-f_I.
\end{cases}
\end{gather*}
Here we can appreciate how this system becomes easier if we add a real coefficient to the kinetics term (as it is the case of Ginzburg-Landau equations) since then it appears a new term $\frac{\partial u_R}{\partial t}$ in the first equation and a new term $\frac{\partial u_I}{\partial t}$ in the second equation. See the paper \cite{MR3208711}.
\end{rmk}

\section{Existence and uniqueness of the solutions}
\label{seceu}

\begin{defi}
\label{defsol}
Let $\Omega\subseteq\R^N$ be an open subset, $(a,b)\in\C^2$ and $V\in L^\infty(\Omega).$
\begin{enumerate}
\item
\label{defsol1}
Let $F\in H^{-1}(\Omega)+L^\infty(\Omega).$ We shall say that a function $u$ is a \textit{global weak solution} to~\eqref{nls} with boundary condition~\eqref{dir}, if $u\in H^1_0(\Omega)\cap L^1(\Omega),$ there is saturated section $U$ associated to $u,$ and if
\begin{multline}
\label{defsol11}
\langle\nabla u,\nabla v\rangle_{L^2(\Omega),L^2(\Omega)}+\langle a\,U,v\rangle_{L^\infty(\Omega),L^1(\Omega)}
												+\langle b\,u,v\rangle_{L^2(\Omega),L^2(\Omega)}	\\
+\langle V\,u,v\rangle_{L^2(\Omega),L^2(\Omega)}=\langle F,v\rangle_{X^\star,X},
\end{multline}
for any $v\in H^1_0(\Omega)\cap L^1(\Omega),$ where $X=H^1_0(\Omega)\cap L^1(\Omega).$
\item
\label{defsol2}
Assume that $\Omega$ has a finite measure and a Lipschitz continuous boundary. Let $F\in H^1(\Omega)^\star.$ We shall say that a function $u$ is a \textit{global weak solution} to \eqref{nls} with boundary condition \eqref{neu} if $u\in H^1(\Omega),$ there is a saturated section $U$ associated to $u,$ and if $(u,U)$ satisfies \eqref{defsol11} for any $v\in H^1(\Omega),$ where $X=H^1(\Omega).$
\end{enumerate}
Sometimes, we shall write $(u,U)$ to designate a solution with the obvious meanings.
\end{defi}

\noindent
By convention, throughout this paper $\Omega$ denotes any open subset of $\R^N,$ and $(a,b)$ is a pair of complex numbers. When a function will be said to satisfy the boundary condition \eqref{neu}, it will always be assumed that $\Omega$ has a finite measure and a Lipschitz continuous boundary.

\noindent
Let,
\begin{gather}
\B=\C\setminus\left\{z\in\C; \Re(z)\le-\frac1{\CP^2} \text{ and } \Im(z)=0\right\},
\end{gather}
where $\CP$ is the constant in Poincaré's inequality~\eqref{poincare} below.

\begin{thm}[\textbf{Existence and a priori bound}]
\label{thmexi1}
Assume that $|\Omega|<\infty$ and $b\in\B.$ Let $V\in L^\infty(\Omega;\R)$ with $V\ge0,$ a.e.\,in $\Omega.$ Then for any $F\in H^{-1}(\Omega),$ equations \eqref{nls}--\eqref{dir} admit at least one global weak solution. In addition, Symmetry Property~$\ref{sym}$ below holds. Finally, any solution $u$ to \eqref{nls}--\eqref{dir} satisfies
\begin{gather}
\label{thmbound1}
\|u\|_{H^1_0(\Omega)}\le C,
\end{gather}
where $C=C(\|F\|_{H^{-1}(\Omega)},\|V\|_{L^\infty(\Omega;\R)},|\Omega|,|a|,|b|,N).$
\end{thm}

\begin{sym}
\label{sym}
Furthermore, if there exists $\vR\in SO_N(\R)$ such that for almost every $x\in\Omega,$ $\vR x\in\Omega,$ $F(\vR x)=F(x)$ and $V(\vR x)=V(x)$ then we may construct a solution $u$ which also satisfies $u(\vR x)=u(x),$ for almost every $x\in\Omega.$ When $N=1,$ if $\Omega$ is symmetric with respect to the origin and if $F$ and $V$ are odd functions then $u$ is also an odd function.
\end{sym}

\noindent
Here and in what follows, $SO_N(\R)$ denotes the special orthogonal group of $\R^N.$ We recall that $\A$ is defined by \eqref{A}.

\begin{thm}[\textbf{Existence and a priori bound}]
\label{thmexi2}
Let $V\in L^\infty(\Omega).$ Assume that $a\in\A,$ $\Im(b)\neq0,$ $\Im(a)\Im(b)\ge0$ and $\Im(b)\Im(V)\ge0,$ a.e.\,in $\Omega.$ Then for any $F\in H^{-1}(\Omega),$ equations~\eqref{nls}--\eqref{dir} admit at least one global weak solution. In addition, Symmetry Property~$\ref{sym}$ holds. Finally, any solution $u$ to \eqref{nls}--\eqref{dir} satisfies
\begin{gather}
\label{thmbound2}
\|u\|_{H^1_0(\Omega)}^2+\|u\|_{L^1(\Omega)}+\int_\Omega|\Im(V)||u|^2\d x\le C\|F\|_{H^{-1}(\Omega)}^2,
\end{gather}
where $C=C(\|\Re(V)\|_{L^\infty(\Omega)},|a|,|b|).$ When $F\in H^1(\Omega)^\star,$ a similar statement holds for the boundary condition \eqref{neu}.
\end{thm}

\begin{rmk}
\label{rmkthmexi2}
Note that if, in addition, $\Re(a)\ge0$ and $\Re(a\ovl b)+\Re(a\ovl V)\ge0,$ a.e.\,in $\Omega,$ then the solution given by Theorem~\ref{thmexi2} is unique (\cite[Theorem~2.8]{MR4848642}).
\end{rmk}

\begin{thm}[\textbf{Null solution}]
\label{thmunull}
Let $V\in L^\infty(\Omega).$ Assume that $a\in\A,$ $\Im(b)\neq0,$ $\Im(a)\Im(b)\ge0$ and $\Im(b)\Im(V)\ge0,$ a.e.\,in $\Omega.$ Then there exists $M=M(|a|,|b|,\|\Re(V)\|_{L^\infty(\Omega)})$ satisfying the following property. Let $F\in L^\infty(\Omega)$ with $\|F\|_{L^\infty(\Omega)}\le|a|.$ If $\|F\|_{L^\infty(\Omega)}\le\frac1M$ then the unique global weak solution $(u,U)$ to~\eqref{nls} with boundary condition~\eqref{dir} or \eqref{neu} is given by,
\begin{gather}
\label{thmunull1}
u=0 \; \text{ and } \; U=\frac1aF,
\end{gather}
almost everywhere in $\Omega.$
\end{thm}

\section{Setting of the framework and proofs of the existence theorems}
\label{secset}

Let $\delta\in\{0,1\}$ and $V\in L^\infty(\Omega).$ For $n\in\N$ and $u\in L^2(\Omega),$ let
\begin{gather}
\label{gn}
g_n(u)=
\begin{cases}
\dfrac{u}{|u|+(n-|u|)\frac1{n^2}},		&	\text{if } |u|\le n, \medskip \\
\dfrac{u}{|u|},					&	\text{if } |u| > n,
\end{cases}
\end{gather}
\begin{gather}
\label{hn}
h_n(u)=
\begin{cases}
u,				&	\text{if } |u|\le n, \medskip \\
n\dfrac{u}{|u|},		&	\text{if } |u| > n,
\end{cases}
\end{gather}
\begin{gather}
\label{fn}
f_{n,\delta}=ag_n(u)+(b-\delta+V)h_n(u).
\end{gather}
Let $X=H^1_0(\Omega)$ if we deal with the boundary condition \eqref{dir}, and let $X=H^1(\Omega)$ if we deal with the boundary condition \eqref{neu}.
\\
Let $F\in X^\star.$ Throughout this section, $u$ denotes any global weak solution to
\begin{gather}
\label{nls1}
-\Delta u+a\,U+b\,u+V\,u=F,
\end{gather}
with boundary condition \eqref{dir} or \eqref{neu}. Moreover, for each $n\in\N,$ $u_n\in H^1_0(\Omega)$ denotes any global weak solution to
\begin{gather}
\label{nls20}
-\Delta u_n+f_{n,0}(u_n)=F,
\end{gather}
with boundary condition \eqref{dir}, and $v_n$ denotes any global weak solution to
\begin{gather}
\label{nls21}
-\Delta v_n+v_n+f_{n,1}(v_n)=F,
\end{gather}
with boundary condition \eqref{dir} or \eqref{neu}. Choosing as test functions $u$ and $\vi u$ in \eqref{nls1}, $u_n$ and $\vi u_n$ in \eqref{nls20}, and $v_n$ and $\vi v_n$ in \eqref{nls21}, we obtain
\begin{gather}
\label{est1}
\|\nabla u\|_{L^2(\Omega)}^2+\Re(a)\|u\|_{L^1(\Omega)}+\Re(b)\|u\|_{L^2(\Omega)}^2+\vint_\Omega\Re(V)|u|^2\d x=\langle F,u\rangle_{X^\star,X},	\\
\label{est2}
\Im(a)\|u\|_{L^1(\Omega)}+\Im(b)\|u\|_{L^2(\Omega)}^2+\vint_\Omega\Im(V)|u|^2\d x=\langle F,\vi u\rangle_{X^\star,X},
\end{gather}
and for any $n\in\N,$
\begin{multline}
\label{estn1}
\|\nabla u_n\|_{L^2(\Omega)}^2
+\Re(a)\left(\int_{\{|u_n|\le n\}}\dfrac{|u_n|^2}{|u_n|+(n-|u_n|)\frac1{n^2}}\d x+\|u_n\|_{L^1(\{|u_n|>n\})}\right)	\\
+\Re(b)\left(\|u_n\|_{L^2(\{|u_n|\le n\})}^2+n\|u_n\|_{L^1(\{|u_n|>n\})}\right)								\\
+\int_{\{|u_n|\le n\}}\Re(V)|u_n|^2\d x+n\int_{\{|u_n|>n\}}\Re(V)|u_n|\d x=\langle F,u_n\rangle_{X^\star,X},
\end{multline}
\begin{multline}
\label{estn2}
\Im(a)\left(\int_{\{|u_n|\le n\}}\dfrac{|u_n|^2}{|u_n|+(n-|u_n|)\frac1{n^2}}\d x+\|u_n\|_{L^1(\{|u_n|>n\})}\right)		\\
+\Im(b)\left(\|u_n\|_{L^2(\{|u_n|\le n\})}^2+n\|u_n\|_{L^1(\{|u_n|>n\})}\right)								\\
+\int_{\{|u_n|\le n\}}\Im(V)|u_n|^2\d x+n\int_{\{|u_n|>n\}}\Im(V)|u_n|\d x=\langle F,\vi u_n\rangle_{X^\star,X},
\end{multline}
\begin{multline}
\label{estn11}
\|v_n\|_X^2
+\Re(a)\left(\int_{\{|v_n|\le n\}}\dfrac{|v_n|^2}{|v_n|+(n-|v_n|)\frac1{n^2}}\d x+\|v_n\|_{L^1(\{|v_n|>n\})}\right)		\\
\le\big(|\Re(b)|+1+\|\Re(V)\|_{L^\infty(\Omega)}\big)\left(\|v_n\|_{L^2(\{|v_n|\le n\})}^2+n\|v_n\|_{L^1(\{|v_n|>n\})}\right)+\langle F,v_n\rangle_{X^\star,X},\end{multline}
and
\begin{gather}
\label{estn21}
v_n \text{ satisfies } \eqref{estn2}.
\end{gather}
We note that for any $w\in L^2(\Omega)$ and $n\in\N,$ we have that
\begin{gather}
\label{est3}
\int_{\{|w|\le n\}}\dfrac{|w|^2}{|w|+(n-|w|)\frac1{n^2}}\d x+\|w\|_{L^1(\{|w|>n\})}\le\|w\|_{L^1(\Omega)}, \\
\label{est4}
\|w\|_{L^2(\{|w|\le n\})}^2+n\|w\|_{L^1(\{|w|>n\})}\le\|w\|_{L^2(\Omega)}^2.
\end{gather}
Finally, we recall that if $|\Omega|<\infty$ then we have Poincaré's inequality:
\begin{gather}
\label{poincare}
\forall w\in H^1_0(\Omega), \; \|w\|_{L^2(\Omega)}\le\CP\|\nabla w\|_{L^2(\Omega)},
\end{gather}
where $\CP=\CP(|\Omega|,N),$ and then
\begin{gather}
\label{L1grad}
\forall w\in H^1_0(\Omega), \; \|w\|_{L^1(\Omega)}\le|\Omega|^\frac12\|w\|_{L^2(\Omega)}\le\CP|\Omega|^\frac12\|\nabla w\|_{L^2(\Omega)},	\\
\label{h10}
\forall w\in H^1_0(\Omega), \; \|w\|_{H^1_0(\Omega)}\le(1+\CP)\|\nabla w\|_{L^2(\Omega)}.
\end{gather}

\subsection{Homogeneous Dirichlet boundary condition with a domain of finite measure}
\label{subsecfin}

Throughout this subsection, we deal with the boundary condition \eqref{dir} and assume that $|\Omega|<\infty.$

\setcounter{equation}{0}
\numberwithin{equation}{subsection}

\begin{vlem}
\label{lemest1}
If $\Re(b)\ge0$ and $\Re(V)\ge0$ then
\begin{gather}
\label{lemest11}
\|\nabla u\|_{L^2(\Omega)}+\|\nabla u_n\|_{L^2(\Omega)}+\vint_\Omega\Re(V)|u|^2\d x\le C(\|F\|_{H^{-1}(\Omega)},|\Omega|,|\Re(a)|,N),
\end{gather}
for any $n\in\N.$
\end{vlem}

\begin{proof*}
Starting with \eqref{estn1} and using \eqref{est3}--\eqref{h10}, we get for any $n\in\N,$
\begin{gather*}
\|\nabla u_n\|_{L^2(\Omega)}^2\le\left(|\Re(a)|\CP|\Omega|^\frac12+(1+\CP)\|F\|_{H^{-1}(\Omega)}\right)\|\nabla u_n\|_{L^2(\Omega)},
\end{gather*}
from which the result follows for $\|\nabla u_n\|_{L^2(\Omega)}.$ Starting with \eqref{est1}, we get the estimate for $\|\nabla u\|_{L^2(\Omega)}+\int_\Omega\Re(V)|u|^2\d x$ in the same way.
\medskip
\end{proof*}

\begin{vlem}
\label{lemest2}
If $-\frac1{\CP^2}<\Re(b)<0$ and $\Re(V)\ge0$ then $u$ and $(u_n)_{n\in\N}$ satisfy \eqref{lemest11}.
\end{vlem}

\begin{proof*}
Starting with \eqref{estn1} and using \eqref{est3}, \eqref{est4}, \eqref{L1grad} and \eqref{h10}, we get for any $n\in\N,$
\begin{gather*}
\|\nabla u_n\|_{L^2(\Omega)}^2\le(|\Re(a)|\CP|\Omega|^\frac12+(1+\CP)\|F\|_{H^{-1}(\Omega)})\|\nabla u_n\|_{L^2(\Omega)}-\Re(b)\CP^2\|\nabla u_n\|_{L^2(\Omega)}^2,
\end{gather*}
from which we get,
\begin{gather*}
(1+\Re(b)\CP^2)\|\nabla u_n\|_{L^2(\Omega)}\le(|\Re(a)|\CP|\Omega|^\frac12+(1+\CP)\|F\|_{H^{-1}(\Omega)}),
\end{gather*}
for any $n\in\N.$ But $1+\Re(b)\CP^2>0$ and then the result follows for $\|\nabla u_n\|_{L^2(\Omega)}.$ Starting with \eqref{est1}, we get the estimate for $\|\nabla u\|_{L^2(\Omega)}+\int_\Omega\Re(V)|u|^2\d x$ in the same way.
\medskip
\end{proof*}

\begin{vlem}
\label{lemest3}
If $\Im(b)\neq0$ and $\Im(b)\Im(V)\ge0$ then for any $n\in\N,$
\begin{gather*}
\|\nabla u\|_{L^2(\Omega)}+\|\nabla u_n\|_{L^2(\Omega)}+\vint_\Omega|\Im(V)||u|^2\d x\le C,
\end{gather*}
where $C=C(\|F\|_{H^{-1}(\Omega)},\|\Re(V)\|_{L^\infty(\Omega)},|\Omega|,|a|,|b|,N).$
\end{vlem}

\begin{proof*}
Let $n\in\N.$ Since $\Im(b)\neq0$ and $\Im(b)\Im(V)\ge0,$ we infer from \eqref{estn2} with help of \eqref{est3}, \eqref{L1grad} and \eqref{h10} that
\begin{gather}
\label{demlemest31}
|\Im(b)|\left(\|u_n\|_{L^2(\{|u_n|\le n\})}^2+n\|u_n\|_{L^1(\{|u_n|>n\})}\right)\le C_1\|\nabla u_n\|_{L^2(\Omega)},
\end{gather}
where
\begin{gather*}
C_1=|\Im(a)||\Omega|^\frac12\CP+(1+\CP)\|F\|_{H^{-1}(\Omega)}.
\end{gather*}
It follows that,
\begin{align*}
	&	\; \int_{\{|u_n|\le n\}}|\Re(V)||u_n|^2\d x+n\int_{\{|u_n|>n\}}|\Re(V)||u_n|\d x						\\
  \le	&	\; \|\Re(V)\|_{L^\infty(\Omega)}\left(\|u_n\|_{L^2(\{|u_n|\le n\})}^2+n\|u_n\|_{L^1(\{|u_n|>n\})}\right)	\\
  \le	&	\; C_1|\Im(b)|^{-1}\|\Re(V)\|_{L^\infty(\Omega)}\|\nabla u_n\|_{L^2(\Omega)}.
\end{align*}
This yields with \eqref{estn1}, \eqref{est3}, \eqref{L1grad}, \eqref{h10} and \eqref{demlemest31} that
\begin{multline*}
\|\nabla u_n\|_{L^2(\Omega)}^2	\le\Big(\CP|\Re(a)||\Omega|^\frac12+C_1|\Im(b)|^{-1}(|\Re(b)|+\|\Re(V)\|_{L^\infty(\Omega)})	\\
+(1+\CP)\|F\|_{H^{-1}(\Omega)}\Big)\|\nabla u_n\|_{L^2(\Omega)},
\end{multline*}
which gives the desired result for $\|\nabla u_n\|_{L^2(\Omega)}.$ For $\|\nabla u\|_{L^2(\Omega)},$ we proceed as follows. Using \eqref{est2} in place of \eqref{est3}, we obtain in the same way as for \eqref{demlemest31} that
\begin{gather}
\label{demlemest32}
|\Im(b)|\|u\|_{L^2(\Omega)}^2+\vint_\Omega|\Im(V)||u|^2\d x\le C_2\|\nabla u\|_{L^2(\Omega)},
\end{gather}
where
\begin{gather*}
C_2=|\Im(a)||\Omega|^\frac12\CP+(1+\CP)\|F\|_{H^{-1}(\Omega)}.
\end{gather*}
As a consequence,
\begin{gather}
\label{demlemest33}
\int_\Omega|\Re(V)||u|^2\d x\le C_2|\Im(b)|^{-1}\|\Re(V)\|_{L^\infty(\Omega)}\|\nabla u\|_{L^2(\Omega)}.
\end{gather}
Using \eqref{L1grad}, \eqref{h10}, \eqref{demlemest32} and \eqref{demlemest33} in \eqref{est1}, we get that
\begin{multline*}
\|\nabla u\|_{L^2(\Omega)}^2\le\Big(\CP|\Re(a)||\Omega|^\frac12+C_2|\Im(b)|^{-1}(|\Re(b)|+\|\Re(V)\|_{L^\infty(\Omega)})	\\
+(1+\CP)\|F\|_{H^{-1}(\Omega)}\Big)\|\nabla u\|_{L^2(\Omega)}.
\end{multline*}
Hence the result with help of \eqref{demlemest32}.
\medskip
\end{proof*}

\subsection{General case}
\label{subsecinf}

In this subsection, we deal with both boundary conditions \eqref{dir} and \eqref{neu}. In addition, no assumption about the open set $\Omega$ is made. We recall that $X=H^1_0(\Omega)$ if we deal with the boundary condition \eqref{dir}, and $X=H^1(\Omega)$ if we deal with the boundary condition \eqref{neu}. Let $F\in X^\star.$

\begin{vlem}
\label{lemest4}
If $a\in\A,$ $\Im(b)\neq0,$ $\Im(a)\Im(b)\ge0$ and $\Im(b)\Im(V)\ge0$ then
\begin{gather}
\label{lemest41}
\|u\|_X^2+\|u\|_{L^1(\Omega)}+\int_\Omega|\Im(V)||u|^2\d x\le C(|\langle F,\vi u\rangle_{X^\star,X}|+|\langle F,u\rangle_{X^\star,X}|),	\\
\label{lemest42}
\|v_n\|_X^2+\int_{\{|v_n|\le n\}}\dfrac{|v_n|^2}{|v_n|+(n-|v_n|)\frac1{n^2}}\d x\le C\|F\|_{X^\star}^2,
\end{gather}
for any $n\in\N,$ where $C=C(\|\Re(V)\|_{L^\infty(\Omega)},|a|,|b|).$
\end{vlem}

\begin{proof*}
Let $n\in\N.$ By our assumptions, \eqref{estn21} may be written as,
\begin{multline}
\label{demlemest41}
|\Im(a)|\left(\int_{\{|v_n|\le n\}}\dfrac{|v_n|^2}{|u_n|+(n-|v_n|)\frac1{n^2}}\d x+\|v_n\|_{L^1(\{|v_n|>n\})}\right)		\\
+|\Im(b)|\left(\|v_n\|_{L^2(\{|v_n|\le n\})}^2+n\|v_n\|_{L^1(\{|v_n|>n\})}\right)								\\
+\int_{\{|v_n|\le n\}}|\Im(V)||u_n|^2\d x+n\int_{\{|v_n|>n\}}|\Im(V)||v_n|\d x=|\langle F,\vi v_n\rangle_{X^\star,X}|,
\end{multline}
If $\Re(a)>0$ then by \eqref{estn11} and \eqref{demlemest41}, we have
\begin{align*}
	&	\; \|v_n\|_X^2+\Re(a)\int_{\{|v_n|\le n\}}\dfrac{|v_n|^2}{|v_n|+(n-|v_n|)\frac1{n^2}}\d x		\\
  \le	&	\; \left(\frac{|\Re(b)|+1+\|\Re(V)\|_{L^\infty(\Omega)}}{|\Im(b)|}\right)|\langle F,\vi v_n\rangle_{X^\star,X}|+|\langle F,v_n\rangle_{X^\star,X}|.
\end{align*}
If $\Re(a)\le0$ then $\Im(a)\neq0.$ Multiplying \eqref{demlemest41} by $L\eqdef\frac{|\Re(a)|+1}{|\Im(a)|}$ and adding the result to \eqref{estn11}, we get that
\begin{align*}
	&	\; \|v_n\|_X^2+\int_{\{|v_n|\le n\}}\dfrac{|v_n|^2}{|v_n|+(n-|v_n|)\frac1{n^2}}\d x		\\
  \le	&	\; \left(\frac{|\Re(b)|+1+\|\Re(V)\|_{L^\infty(\Omega)}}{|\Im(b)|}+L\right)|\langle F,\vi v_n\rangle_{X^\star,X}|+|\langle F,v_n\rangle_{X^\star,X}|.
\end{align*}
In both cases we obtain that,
\begin{gather*}
\|v_n\|_X^2+\int_{\{|v_n|\le n\}}\dfrac{|v_n|^2}{|u_n|+(n-|v_n|)\frac1{n^2}}\d x\le C(|\langle F,\vi v_n\rangle_{X^\star,X}|+|\langle F,v_n\rangle_{X^\star,X}|),
\end{gather*}
for some $C=C(\|\Re(V)\|_{L^\infty(\Omega)},|a|,|b|).$ Applying Young's inequality to the above, we get \eqref{lemest42}. Using \eqref{est1} and \eqref{est2} instead of \eqref{estn11} and \eqref{estn21}, we obtain \eqref{lemest41} in the same way.
\medskip
\end{proof*}

\subsection{Proofs of the existence and compactness theorems}
\label{subsecproof}

\begin{vproof}{of Theorems~\ref{thmexi1} and \ref{thmexi2}.}
We first note that \eqref{thmbound1} comes from Lemmas~\ref{lemest1}--\ref{lemest3} and \eqref{poincare}, and that \eqref{thmbound2} comes from Lemma~\ref{lemest4} and Young's inequality. It remains to establish the existence part of the theorems. We first assume that $|\Omega|<\infty.$ Let $F$ be as in the theorems. For each $n\in\N,$ let $u_n$ be a global weak solution to \eqref{nls20} and \eqref{dir}, and let $v_n$ be a global weak solution to \eqref{nls21} and \eqref{dir} (respectively, to \eqref{nls21} and \eqref{neu}). Indeed, such solutions exist with help of \cite[Lemma~6.5]{MR4848642}. By Lemmas~\ref{lemest1}--\ref{lemest3}, \eqref{est3}, \eqref{poincare} and \eqref{L1grad}, it follows that $(u_n)_{n\in\N}$ is bounded in $H^1_0(\Omega)$ and
\begin{gather*}
\left(\dfrac{|u_n|^2}{|u_n|+(n-|u_n|)\frac1{n^2}}\1_{\{|u_n|\le n\}}\right)_{n\in\N} \; \text{ is bounded in } L^1(\Omega).
\end{gather*}
By \cite[Lemma~6.2]{MR4848642}, we may extract a subsequence of $(u_n)_{n\in\N}$ which converges to a solution of \eqref{nls}--\eqref{dir}. Theorem~\ref{thmexi1} is then proved. By Lemma~\ref{lemest4}, $(v_n)_{n\in\N}$ is bounded in $H^1_0(\Omega)$ (respectively, in $H^1(\Omega))$ and
\begin{gather*}
\left(\dfrac{|v_n|^2}{|v_n|+(n-|v_n|)\frac1{n^2}}\1_{\{|v_n|\le n\}}\right)_{n\in\N} \; \text{ is bounded in } L^1(\Omega).
\end{gather*}
By \cite[Lemma~6.2]{MR4848642}, (respectively, \cite[Lemma~6.3]{MR4848642},) we may extract a subsequence of $(v_n)_{n\in\N}$ which converges to a solution of \eqref{nls}--\eqref{dir} (respectively, of \eqref{nls} and \eqref{neu}). This completes the proof of Theorem~\ref{thmexi1}, then Theorem~\ref{thmexi2} is then proved in the case $|\Omega|<\infty.$ To complete the proof, it remains to show that \eqref{nls}--\eqref{dir} admits a solution when $|\Omega|=\infty.$ An appeal to \eqref{thmbound2} and the Extension Lemma (\cite[Lemma~6.9]{MR4848642} applied with $\Omega_n=\Omega\cap B(0,n))$ gives the existence of a $u\in H^1_0(\Omega)$ and of a saturated section $U$ associated to $u$ such that $(u,U)$ satisfies \eqref{nls} in $\Dr^\p(\Omega).$ But $\Delta u,Vu,F\in H^{-1}(\Omega)$ and $U\in L^\infty(\Omega)$ so that the equation \eqref{nls} makes sense in $H^{-1}(\Omega)+L^\infty(\Omega)\inj\Dr^\p(\Omega).$ Theorem~\ref{thmexi2} is then proved.
\medskip
\end{vproof}

\begin{vproof}{of Theorem~\ref{thmunull}.}
We indeed check that $(u,U)$ defined by \eqref{thmunull1} is a solution to \eqref{nls}. Now, assume that $(u,U)$ is a solution to \eqref{nls}. Taking the duality product of \eqref{nls} with $u$ and $\vi u,$ we have that,
\begin{gather}
\label{demthmunull1}
\|\nabla u\|_{L^2(\Omega)}^2+\Re(a)\|u\|_{L^1(\Omega)}+(\Re(b)-\|V\|_{L^\infty(\Omega)})\|u\|_{L^2(\Omega)}^2\le\int_\Omega|Fu|\d x,	\\
\label{demthmunull2}
|\Im(a)|\|u\|_{L^1(\Omega)}+|\Im(b)|\|u\|_{L^2(\Omega)}^2+\int_\Omega|\Im(V)||u|^2\d x\le\vint_\Omega|Fu|\d x.
\end{gather}
Since we have either $\Re(a)>0$ or $|\Im(a)|>0,$ and since $|\Im(b)|>0,$ we may find a $C=C(|a|,|b|,\|\Re(V)\|_{L^\infty(\Omega)})$ such that $\Re(a)+C|\Im(a)|>0$ and $\Re(b)-\|V\|_{L^\infty(\Omega)}+C|\Im(b)|\ge1.$ We then multiply \eqref{demthmunull2} by $C$ and sum the result to \eqref{demthmunull1}. This yields to,
\begin{gather*}
\|u\|_{H^1(\Omega)}^2+\|u\|_{L^1(\Omega)}\le M\vint_\Omega|Fu|\d x,
\end{gather*}
for some $M=M(|a|,|b|,\|\Re(V)\|_{L^\infty(\Omega)}).$ Applying Hölder's inequality to the above, we get that
\begin{gather*}
\|u\|_{H^1(\Omega)}^2+(1-M\|F\|_{L^\infty(\Omega)})\|u\|_{L^1(\Omega)}\le0.
\end{gather*}
Hence \eqref{thmunull1} if $\|F\|_{L^\infty(\Omega)}\le\frac1M.$
\medskip
\end{vproof}

\begin{vproof}{of Theorem~\ref{thmmain}.}
Let $K$ be a compact subset of $\R^N$ for which $F_{|K^\co}\in L^\infty(K^\co).$ Let $R>0$ be such that $K\subset\ovl B(0,R)$ and let $\eps\in(0,1).$ \\
\textbf{Proof of Property~\ref{thmmain1}.}
Let us write \eqref{g} as follows.
\begin{gather}
\label{demthmmain1}
-\Delta g+a\,G+bg+Vg=F_1,
\end{gather}
where $b=-\vi\frac{N+2p}4,$ $V(x)=-\frac1{16}|x|^2$ and $F_1=-Fe^{-\vi\frac{|x|^2}8}.$ We have that $\Im(b)=-\frac{N+4}4<0,$ $\Im(a)\Im(b)\ge0$ and $\Im(b)\Im(V)=0,$ in $\R^N.$ It follows that \eqref{demthmmain1} falls into the scope of Theorem~\ref{thmexi2} and then \eqref{demthmmain1} admits a solution $g_\eps\in H^1_0(B(0,R+2\eps)),$ where the right member of \eqref{demthmmain1} is $F_{1|B(0,R+2\eps)}.$ By global elliptic regularity $g_\eps\in H^2(B(0,R+2\eps))$ (Gilbarg and Trudinger~\cite[Theorem~8.12, p.186]{MR1814364}). Let us denote by $G_\eps$ the saturated section associated to $g_\eps.$ Applying \cite[Theorem~3.1]{MR3190983}, we have that,
\begin{multline}
\label{demthmmain2}
\|\nabla g_\eps\|_{L^2(B_{R,\eps,x_0}(\rho))}^2+\Re(a)\|g_\eps\|_{L^1(B_{R,\eps,x_0}(\rho))}
								+\Re(b)\|g_\eps\|_{L^2(B_{R,\eps,x_0}(\rho))}^2		\\
-\vint_{B_{R,\eps,x_0}(\rho)}\frac{|x|^2}{16}|g_\eps|^2\d x
=\Re\left(\:\vint_{B_{R,\eps,x_0}(\rho)}F_1\,\ovl g_\eps\,\d x\right)+\Re\left(\:\vint_{\vsS_{R,\eps,x_0}(\rho)}g_\eps\ovl{\nabla g_\eps}.\frac{x-x_0}{|x-x_0|}\d\sigma\right),
\end{multline}
\begin{multline}
\label{demthmmain3}
|\Im(a)|\|g_\eps\|_{L^1(B_{R,\eps,x_0}(\rho))}+|\Im(b)|\|g_\eps\|_{L^2(B_{R,\eps,x_0}(\rho))}^2		\\
=-\Im\left(\:\vint_{B_{R,\eps,x_0}(\rho)}F_1\,\ovl g_\eps\,\d x\right)-\Im\left(\:\vint_{\vsS_{R,\eps,x_0}(\rho)}g_\eps\ovl{\nabla g_\eps}.\frac{x-x_0}{|x-x_0|}\d\sigma\right),
\end{multline}
for any $x_0\in B(0,R+2\eps)$ and $\rho\in[0,2\eps),$ where $B_{R,\eps,x_0}(\rho)=B(0,R+2\eps)\cap B(x_0,\rho)$ and $\vsS_{R,\eps,x_0}(\rho)=B(0,R+2\eps)\cap\vsS(x_0,\rho).$ Let us denote by $g\in H^1(\R^N)$ the extension by $0$ of $g_\eps$ outside of $B(0,R+2\eps).$ Since we have either $\Re(a)>0$ or $|\Im(a)|>0,$ and $|\Im(b)|>0,$ we may find a $C=C(|a|,|\Im(p)|,R,N)$ such that $\Re(a)+C|\Im(a)|>0$ and $\Re(b)-\frac{(R+2)^2}{16}+C|\Im(b)|\ge1.$ We then multiply \eqref{demthmmain3} by $C$ and sum the result to \eqref{demthmmain2}. This yields to,
\begin{gather*}
\|g\|_{H^1(B(x_0,\rho))}^2+\|g\|_{L^1(B(x_0,\rho))}
\le C_1\left(\dsp\int_{B(x_0,\rho)}|F_1g|+\left|\dsp\int_{\vsS(x_0,\rho)}g\ovl{\nabla g}.\frac{x-x_0}{|x-x_0|}\d\sigma\right|\right),
\end{gather*}
for any $x_0\in B(0,R+2\eps)$ and $\rho\in[0,2\eps),$ and for some $C_1=C_1(|a|,|\Im(p)|,R,N).$ It follows from Hölder's inequality that for $M\ge2C_1,$ if $\|F\|_{L^\infty(K^\co)}\le\frac1M$ then
\begin{gather}
\label{demthmmain4}
\|g\|_{H^1(B(x_0,\rho))}^2+\|g\|_{L^1(B(x_0,\rho))}
\le M\left|\dsp\int_{\vsS(x_0,\rho)}g\ovl{\nabla g}.\frac{x-x_0}{|x-x_0|}\d\sigma\right|,
\end{gather}
for any $x_0\in B(0,R+2\eps)$ and $\rho\in[0,2\eps)$ such that $K\cap B(x_0,2\eps)=\emptyset.$ It follows from \cite[Theorem~4.1]{MR4848642} that there exists $\rho_\M\ge0$ such that $g=0,$ a.e.\,in $B(x_0,\rho_\M),$ for any $x_0\in B(0,R+2\eps)$ such that $K\cap B(x_0,2\eps)=\emptyset.$ By \eqref{thmbound2} and \cite[Theorem~4.1]{MR4848642}, there exists $\delta=\delta(|a|,|\Im(p)|,R,\eps,N)$ such that if $\|F\|_{L^2(\R^N)}\le\delta$ then $\rho_\M>\eps.$ We then deduce that $g=g_\eps=0,$ a.e.\,in $B(0,R+2\eps)\setminus K(\eps).$ Now, let us define $G$ on $\R^N$ by $G=G_\eps,$ in $B(0,R+2\eps)$ and by $G=-\frac1aFe^{-\vi\frac{|x|^2}8},$ in $B(0,R+2\eps)^\co.$ Choosing also $M\ge|a|^{-1},$ it follows that $G$ is a saturated section associated to $g.$ So, we have shown that $(g,G)$ is a solution to \eqref{g}, $g\in H^2(\R^N)$ and $\supp g\subset K(\eps).$ Now, we define $\vphi$ and $\Phi$ by \eqref{gU} and \eqref{gP}, respectively, and finally, $u$ and $U$ by \eqref{eqprou} and \eqref{eqproU}, respectively. The proof of \eqref{thmmain11} comes from standards arguments of integration theory but for convenience of the reader, we postpone its proof to the Appendix~\ref{appendix}. This completes the proof.
\\
\textbf{Proof of Property~\ref{thmmain2}.}
Using the change of functions \eqref{gU} and \eqref{gP}, we are brought back to show the uniqueness for the equation \eqref{g}. In both cases \eqref{thmmain21} and \eqref{thmmain22}, $\vphi$ and $\phi$ belong to $L^2(\R^N)$ and are compactly supported. It follows that the corresponding solutions to \eqref{g} belong to $L^2(\R^N)$ and their Laplacian belong to $L^2_\loc(\R^N).$ By interior elliptic regularity, they belong to $H^2_\loc(\R^N)$ (Cazenave~\cite[Proposition~4.1.2]{caz-sle}). Since they are compactly supported, they actually belong to $H^2(\R^N)$ and it is sufficient to show the uniqueness for \eqref{g} set in $B(0,r),$ where $r>0$ is large enough to have $\supp\vphi\cup\supp\phi\subset B(0,r).$ It follows that \eqref{g} falls into the scope of the uniqueness~\cite[Theorem~2.8]{MR4848642}. Since also $a\in\A$ and $\Re(a)\ge0,$ we only have to show that,
\begin{gather*}
\Re(a\ovl b)+\Re(a\ovl V)>0, \; \text{a.e.\,in } B(0,r),
\end{gather*}
where $b$ and $V$ are as in  \eqref{demthmmain1}. If $\Re(a)=0$ then $\Im(a)<0$ and $\Re(a\ovl b)+\Re(a\ovl V)=-\Im(a)\frac{N+4}4>0,$ over $\R^N.$ If $\Re(a)>0$ then
\begin{gather*}
\Re(a\ovl b)+\Re(a\ovl V)=\frac12\Re(a)\Im(p)-\Im(a)\frac{N+4}4-\frac1{16}\Re(a)|x|^2, \; \text{in } \R^N.
\end{gather*}
Using \eqref{thmmain22}, we have that
\begin{gather*}
\Re(a\ovl b)+\Re(a\ovl V)>\frac{\Re(a)}{16}\left(8\Im(p)-4\frac{\Im(a)}{\Re(a)}(N+4)-r^2\right)\ge0, \; \text{in } B(0,r).
\end{gather*}
This concludes the proof of the theorem.
\medskip
\end{vproof}

\appendix
\section{Appendix}
\label{appendix}

\numberwithin{equation}{section}

\begin{lem}
\label{lemA}
Let $m\in\N_0,$ $1<q<\infty$ and $\vphi\in W^{m,q}(\R^N).$ Let $p\in\C,$ and let $u$ be defined by \eqref{eqprou}. Then,
\begin{gather}
\label{lemAA}
u\in C\big((0,\infty);W^{m,q}(\R^N)\big).
\end{gather}
If, in addition, $\supp\vphi$ is compact and $m\ge1$ then,
\begin{gather}
\label{lemAB}
u\in\bigcap_{j=1}^m C^j\big((0,\infty);W^{m-j,q}(\R^N)\big).
\end{gather}
\end{lem}

\begin{proof*}
Let $1<q<\infty,$ $p\in\C,$ $\vphi\in L^q(\R^N)$ and $u$ be defined by \eqref{eqprou}. Let $t>0.$ Let $(t_n)_{n\in\N}\subset(0,\infty)$ be such that $t_n\xrightarrow{n\to\infty}t.$ We claim that,
\begin{gather}
\label{demlemA1}
u(t_n)\underset{n\to\infty}{\overset{L^q(\Omega)_\w}{-\!\!\!-\!\!\!-\!\!\!-\!\!\!-\!\!\!\weak}}u(t).
\end{gather}
By \eqref{rmkthmmain21}, $(u(t_n))_{n\in\N}$ is bounded in $L^q(\R^N).$ So, it is enough to show that $u(t_n)\xrightarrow[n\to\infty]{\Dr^\p(\R^N)}u(t).$ Let $\theta\in\Dr(\R^N).$ By change of variables, we have for any $n\in\N,$
\begin{gather*}
\langle u(t_n),\theta\rangle_{\Dr^\p(\R^N),\Dr(\R^N)}=\Re\int_{\R^N}t_n^\frac{p+N}2\vphi(x)\ovl{\theta(\sqrt{t_n}x)}\d x,	\\
\langle u(t),\theta\rangle_{\Dr^\p(\R^N),\Dr(\R^N)}=\Re\int_{\R^N}t^\frac{p+N}2\vphi(x)\ovl{\theta(\sqrt{t}x)}\d x.
\end{gather*}
It follows from the dominated convergence Theorem that,
\begin{gather*}
\langle u(t_n),\theta\rangle_{\Dr^\p(\R^N),\Dr(\R^N)}\xrightarrow{n\to\infty}\langle u(t),\theta\rangle_{\Dr^\p(\R^N),\Dr(\R^N)}.
\end{gather*}
from which we get \eqref{demlemA1}. By \eqref{rmkthmmain21}, we also have that
\begin{gather}
\label{demlemA2}
\|u(t_n)\|_{L^q(\R^N)}\xrightarrow{n\to\infty}\|u(t)\|_{L^q(\R^N)}.
\end{gather}
By \eqref{demlemA1}, \eqref{demlemA2} and the uniform convexity of the $L^q$-spaces, we infer that
\begin{gather*}
u(t_n)\xrightarrow[n\to\infty]{L^q(\R^N)}u(t),
\end{gather*}
proving that $u\in C\big((0,\infty);L^q(\R^N)\big).$ Now assume that  $\vphi\in W^{m,q}(\R^N),$ for some $m\in\N.$ Then \eqref{lemAA} follows with the same arguments. We have for any $n\in\N$ and almost every $x\in\R^N,$
\begin{gather*}
\frac{\partial u}{\partial t}(t_n,x)=\frac{p}2t_n^\frac{p-2}2\vphi\left(\frac{x}{\sqrt{t_n}}\right)-\frac12t_n^\frac{p-3}{2}x.\nabla\vphi\left(\frac{x}{\sqrt{t_n}}\right),	\\
\frac{\partial u}{\partial t}(t,x)=\frac{p}2t^\frac{p-2}2\vphi\left(\frac{x}{\sqrt{t}}\right)-\frac12t^\frac{p-3}{2}x.\nabla\vphi\left(\frac{x}{\sqrt{t}}\right).
\end{gather*}
If $\supp\vphi$ is compact then we may proceed as above to show that
\begin{gather*}
\frac{p}2t_n^\frac{p-2}2\vphi\left(\frac.{\sqrt{t_n}}\right)\xrightarrow[n\to\infty]{L^q(\R^N)}\frac{p}2t^\frac{p-2}2\vphi\left(\frac{.}{\sqrt{t}}\right),	\\
\frac12t_n^\frac{p-3}{2}(\:.\:).\nabla\vphi\left(\frac{.}{\sqrt{t_n}}\right)\xrightarrow[n\to\infty]{L^q(\R^N)}\frac12t^\frac{p-3}{2}(\:.\:).\nabla\vphi\left(\frac{.}{\sqrt{t}}\right).
\end{gather*}
As a consequence, $\frac{\partial u}{\partial t}(t_n)\xrightarrow[n\to\infty]{L^q(\R^N)}\frac{\partial u}{\partial t}(t)$ and then  $u\in C^1\big((0,\infty);L^q(\R^N)\big).$ The others regularity in \eqref{lemAB} are obtained in the same way and the details are left to the reader.
\medskip
\end{proof*}

\end{document}